\documentclass[a4paper,10pt]{amsart}

\usepackage[cp850]{inputenc}
\usepackage{marvosym}
\usepackage{wasysym}
\usepackage{latexsym}
\usepackage{amsfonts}
\usepackage{mathrsfs}
\usepackage{hyperref}

\newcommand{\R}{\mathbb{R}}
\newcommand{\C}{\mathbb{C}}
\newcommand{\N}{\mathbb{N}}

\newtheorem{theo}{Theorem}
\newtheorem{lem}[theo]{Lemma}

\newtheorem{prop}[theo]{Proposition}

\title[Augmented differential operators]{The augmented operator of a surjective partial differential operator with constant coefficients need not be surjective}
\author{T. Kalmes}

\begin{document}

\begin{abstract}
For $d\geq 3$ we give an example of a constant coefficient surjective differential operator $P(D):\mathscr{D}'(X)\rightarrow\mathscr{D}'(X)$ over some open subset $X\subset\R^d$ such that $P^+(D):\mathscr{D}'(X\times\R)\rightarrow\mathscr{D}'(X\times\R)$ is not surjective, where $P^+(x_1,\ldots,x_{d+1}):=P(x_1,\ldots,x_d)$. This answers in the negative a problem posed by Bonet and Doma\'nski in \cite[Problem 9.1]{Bonet}.
\end{abstract}

\maketitle

\section{Introduction}

For an open subset $X\subset\R^d$ and $P\in\C[X_1,\ldots, X_d]$ a non-zero polynomial consider the corresponding differential operator $P(D)$ on $\mathscr{D}'(X)$, where as usual $D_j=-i\frac{\partial}{\partial x_j}$. For $(x_1,\ldots,x_{d+1})\in\R^{d+1}$ we set $P^+(x_1,\ldots,x_{d+1}):=P(x_1,\ldots,x_d)$ and call $P^+(D)$ the augmented operator, i.e.\ $P(D)$ acting ''on the first $d$ variables'' on $\mathscr{D}'(X\times\R)$.

In \cite[Problem 9.1]{Bonet} Bonet and Doma\'nski asked the natural question whether surjectivity of the constant coefficient differential operator $P(D):\mathscr{D}'(X)\rightarrow\mathscr{D}'(X)$ passes on to surjectivity of $P^+(D):\mathscr{D}'(X\times\R)\rightarrow\mathscr{D}'(X\times\R)$. This problem is loosely connected with the parameter dependence of solutions of the differential equation
\[P(D)u_\lambda=f_\lambda,\]
see \cite{Bonet}. Bonet and Doma\'nski proved in \cite[Proposition 8.3]{Bonet} that for a surjective differential operator $P(D):\mathscr{D}'(X)\rightarrow\mathscr{D}'(X)$ the augmented operator $P^+(D)$ is surjective if and only if the kernel of $P(D)$ has the linear topological invariant $(P\Omega)$. Moreover, a positive answer to this problem would have additional consequences. As shown by Bonet and Doma\'nski in \cite{Bonet} surjectivity of $P^+(D)$ is equivalent to the surjectivity of the vector valued operator $P(D):\mathscr{D}'(X,\Lambda_r(\alpha)')\rightarrow\mathscr{D}'(X,\Lambda_r(\alpha)')$, where $\Lambda_r(\alpha)$ is a power series space, like e.g.\ spaces of smooth functions on a compact manifold, or spaces of holomorphic functions on a nice domain.

Vogt showed in \cite[Proposition 3.4]{Vogt} that the kernel of an elliptic differential operator $P(D)$ always has the property $(\Omega)$. Since for elliptic $P$ the kernels of
\[P(D):C^\infty(X)\rightarrow C^\infty(X)\mbox{ and }P(D):\mathscr{D}'(X)\rightarrow\mathscr{D}'(X)\]
coincide as locally convex spaces it is a Fr\'echet-Schwartz space. Hence it has $(\Omega)$ if and only if it has $(P\Omega)$. So for elliptic operators the problem of Bonet and Doma\'nski always has a positive solution.

By the well-known characterization of surjectivity of constant coefficient differential operators due to H\"ormander, the above problem is equivalent to the question whether $X\times\R$ is $P^+$-convex for supports and singular supports in case of $X$ being $P$-convex for supports and singular supports. In \cite[Proposition 1]{Frerick} it is shown that $P$-convexity for supports of $X$ is passed on to $P^+$-convexity for supports of $X\times\R$. Moreover, it is shown in \cite[Example 9]{Frerick} that an analogous implication for $P$-convexity for singular supports is not true in general but in this example the set $X$ is not $P$-convex for supports.

In \cite{kalmes2} it is shown that in case of $d=2$ the problem has always a positive solution and this is also the case in arbitrary dimension if $P$ is homogeneous, semi-elliptic, or of principal type and $X$ has a special form.\\

The purpose of the present paper is to show that in general the problem posed by Bonet and Doma\'nski has to be answered in the negative. More precisely, for any $d>2$ we present a operator $P(D)$, even hypoelliptic, and an open subset $X\subset\R^d$ such that $P(D)$ is surjective on $\mathscr{D}'(X)$ but $P^+(D)$ is not surjective on $X\times\R$. The polynomial $P$ is inspired by \cite[Example 6.4]{Hoermander singular}.

\section{Notation and auxiliary results}

We begin by introducing some terminology. Recall that a cone $C$ is called {\it proper} if it does not contain any affine subspace of dimension one. Moreover, recall that for an open convex cone $\Gamma\subset\R^d$ its {\it dual cone} is defined as
\[\Gamma^\circ:=\{\xi\in\R^d;\,\forall\,y\in\Gamma:\,\langle y,\xi\rangle\geq 0\}.\]

We will use the following function introduced by H\"ormander in connection with continuation of differentiability (cf.\ \cite[Section 11.3, vol.\ II]{Hoermander}). For a subspace $V$ of $\R^d$ \[\sigma_P(V)=\inf_{t>1}\liminf_{\xi\rightarrow\infty}\tilde{P}_V(\xi,t)/\tilde{P}(\xi,t)\]
with $\tilde{P}_V(\xi,t):=\sup\{|P(\xi+\eta)|;\,\eta\in V,|\eta|\leq t\}, \tilde{P}(\xi,t):=\tilde{P}_{\R^d}(\xi,t)$. This quantity is closely related with the localizations at infinity of the polynomial $P$ which in turn are connected with bounds for the wave front set and the singular support of regular fundamental solutions of $P$. In order to simplify notation we will write $\sigma_P(y)$ instead of $\sigma_P(span\{y\})$. Recall that the localizations of $P$ at infinity are the non-zero multiples of limits of normalized polynomials $x\mapsto P(x+\xi)/\tilde{P}(\xi,1)$ when $\xi$ tends to infinity (see e.g.\ \cite[Section 10.2, vol.\ II]{Hoermander}).

It follows immediately from the definitions that
\begin{equation}\label{localization dominates}
	\sigma_P(V)\leq\inf_{t>1}\frac{\tilde{Q}_V(0,t)}{\tilde{Q}(0,t)}
\end{equation}
for any localization $Q$ of $P$ at infinity and every subspace $V$ of $\R^d$.

Moreover, we will need
\[\sigma^0_P(V):=\inf_{t>1,\xi\in\R^d}\tilde{P}_V(\xi,t)/\tilde{P}(\xi,t).\]
This function has already been considered by H\"ormander in \cite[Section 5]{Hoermander singular} to discuss ``H\"older estimates'' for solutions of partial differential equations. It plays also a crucial r\^ole in our context, as is shown in \cite{Frerick} and \cite{kalmes1}. Obviously, we have
\begin{equation}\label{power}
	\sigma^0_{P^k}(V)=(\sigma^0_P(V))^k
\end{equation}
for every $k\in\N$ as well as
\begin{equation}\label{sigma sigma^0}
	\sigma^0_P(V)\leq\sigma_P(V).
\end{equation}

We will need the following lemma. For a proof, see \cite[Theorem 13]{kalmes2}.

\begin{lem}\label{characterization of p-convexity in the complement of cones}
Let $\Gamma\neq\R^d$ be a non-empty open proper convex cone in $\R^d$ and $X:=\R^d\backslash\Gamma^\circ$. Let $P$ be a non-constant polynomial with principal part $P_m$.
\begin{itemize}
	\item[i)] $X$ is $P$-convex for supports if and only if $P_m(y)\neq 0$ for all $y\in\Gamma$.
	\item[ii)] $X\times\R$ is $P^+$-convex for singular supports if and only if $\sigma^0_P(y)\neq 0$ for all $y\in\Gamma$.
\end{itemize}
\end{lem}

\begin{prop}\label{simple characteristic}
Let $P$ be a homogeneous polynomial of degree $m$ and $\xi\in\R^d$ with $P(\xi)=0$ and $\nabla P(\xi)\neq 0$. Then $x\mapsto\sum_{j=1}^d\partial_j P(\xi)x_j=\langle\nabla P(\xi),x\rangle$ is a localization of $P$ at infinity.
\end{prop}

{\sc Proof.} Taylor expansion about $n\xi$ and the homogeneity of $P$ give for all $x\in\R^d$ and $n\in\N$
\begin{eqnarray*}
	P(x+n\xi)&=&n^{m-1}\left(\langle\nabla P(\xi),x\rangle +\sum_{2\leq|\alpha|\leq m}n^{1-|\alpha|}\frac{P^{(\alpha)}(\xi)}{\alpha !}x^\alpha\right)
\end{eqnarray*}
as well as
\begin{eqnarray*}
	\sqrt{\sum_{0\leq|\alpha|\leq m}|P^{(\alpha)}(n\xi)|^2}&=&n^{m-1}\sqrt{|\nabla P(\xi)|^2+\sum_{2\leq|\alpha|\leq m}n^{2(1-|\alpha|)}|P^{(\alpha)}(\xi)|^2}.
\end{eqnarray*}
Setting
\[a_n:=\sqrt{|\nabla P(\xi)|^2+\sum_{2\leq|\alpha|\leq m}n^{2(1-|\alpha|)}|P^{(\alpha)}(\xi)|^2}\]
and
\[b_n:=\sqrt{1+\sum_{2\leq|\alpha|\leq m}n^{2(1-|\alpha|)}\frac{|P^{(\alpha)}(\xi)|^2}{|\nabla P(\xi)|^2}}\]
it follows
\[\lim_{n\rightarrow\infty}a_n=|\nabla P(\xi)|\mbox{ and }\lim_{n\rightarrow\infty}b_n=1.\]
Moreover, we have
\begin{eqnarray*}
	\frac{\langle\nabla P(\xi),x\rangle}{|\nabla P(\xi)|}&=&\sqrt{\frac{a_n^2}{|\nabla P(\xi)|^2}}\;\frac{\langle\nabla P(\xi),x\rangle}{a_n}\\
	&=&b_n\frac{\langle\nabla P(\xi),x\rangle}{a_n}
\end{eqnarray*}
so that for $x\in\R^d$
\begin{eqnarray*}
	&&|\frac{P(x+n\xi)}{\sqrt{\sum_{0\leq|\alpha|\leq m}|P^{(\alpha)}(n\xi)|^2}}-\frac{\langle\nabla P(\xi),x\rangle}{|\nabla P(\xi)|}|\\
	&=&|\frac{\langle\nabla P(\xi),x\rangle +\sum_{2\leq|\alpha|\leq m}n^{1-|\alpha|}\frac{P^{(\alpha)}(\xi)}{\alpha !}x^\alpha}{a_n}-b_n\frac{\langle\nabla P(\xi),x\rangle}{a_n}|\\
	&\leq& |\frac{1-b_n}{a_n}\;\;\langle \nabla P(\xi),x\rangle|+\sum_{2\leq|\alpha|\leq m}n^{1-|\alpha|}\frac{|P^{(\alpha)}(\xi)|}{\alpha !} \;\;\frac{|x|^m}{a_n}.
\end{eqnarray*}
Therefore,
\[\lim_{n\rightarrow\infty}\sup_{|x|\leq 1}|\frac{P(x+n\xi)}{\sqrt{\sum_{0\leq|\alpha|\leq m}|P^{(\alpha)}(n\xi)|^2}}-\frac{\langle\nabla P(\xi),x\rangle}{|\nabla P(\xi)|}|=0.\]
Applying the equivalence of the norms $Q\mapsto\tilde{Q}(0,1)$ and $Q\mapsto\sqrt{\sum_{0\leq|\alpha|\leq m}|Q^{(\alpha)}(0)|^2}$ on the finite dimensional vector space of polynomials of degree at most $m$ to the polynomials $x\mapsto P(x+n \xi)$ we obtain that for some subsequence $(n_k)_{k\in\N}$ and some $c>0$ we have
\[\lim_{k\rightarrow\infty}\frac{\sqrt{\sum_{0\leq|\alpha|\leq m}|P^{(\alpha)}(n_k\xi)|^2}}{\tilde{P}(n_k\xi,1)}=c.\]
Therefore we obtain
\[\lim_{k\rightarrow\infty}\sup_{|x|\leq 1}|\frac{P(x+n_k\xi)}{\tilde{P}(n_k\xi,1)}-c\frac{\langle\nabla P(\xi),x\rangle}{|\nabla P(\xi)|}|=0\]
so that $x\mapsto \langle\nabla P(\xi),x\rangle$ is a localization of $P$ at infinity.\hfill$\square$\\

The next lemma is \cite[Lemma 6.1]{Hoermander singular}.

\begin{lem}\label{principal dominates}
Let $P$ be a polynomial with principal part $P_m$. Then for any subspace $V\subseteq\R^d$ we have $\sigma^0_P(V)\leq\sigma^0_{P_m}(V)$.
\end{lem}

\section{The example}

As already mentioned in the introduction, the polynomial we are going to construct now is inspired by \cite[Example 6.4]{Hoermander singular}.

Let $Q$ be a homogeneous polynomial of real principal type of degree $m$ and $x\in\R^d,|x|=1$ such that $Q(x)\neq 0$ but $\sigma_Q(x)=0$. By \cite[Lemma 4]{kalmes1} this is only possible when $d\geq 3$. For example, take $Q(\xi)=\xi_1^2-\xi_2^2-\ldots-\xi_d^2$ and $x=e_d=(0,\ldots,0,1)$. Indeed, by Proposition \ref{simple characteristic} applied to $Q$ and $\xi=(1,1,0,\ldots,0)$ it follows that
\[x\mapsto\langle\nabla Q(\xi),x\rangle=2x_1-2x_2\]
is a localization of $Q$ at infinity. Because $\langle\nabla Q(\xi),e_d\rangle=0$ since $d\geq 3$ we have $\sigma_Q(e_d)=0$ by equation (\ref{localization dominates}) and furthermore $Q(e_d)=-1$.

By \cite[Theorem 11.1.12, vol.\ II]{Hoermander} there is a polynomial $R$ of degree $4m-2$ such that
\[P(\xi):=Q(\xi)^{4}+R(\xi)\]
is a hypoelliptic polynomial of degree $4m$. Clearly, for its principal part $P_{4m}$ we have $P_{4m}=Q^{4}$ so $P_{4m}(x)=Q^{4}(x)\neq 0$. On the other hand we also have
\begin{equation}\label{key inequality}
	\sigma^0_P(x)\leq\sigma^0_{P_{4m}}(x)=\sigma^0_{Q^4}(x)=(\sigma^0_Q(x))^4\leq(\sigma_Q(x))^4=0,
\end{equation}
where we have used Lemma \ref{principal dominates}, equation (\ref{power}), and inequality (\ref{sigma sigma^0}).

Because $P_{4m}(x)\neq 0$ and $P_{4m}$ is homogeneous there is an open proper convex cone $\Gamma\neq\R^d$ with $x\in\Gamma$ such that $P_{4m}(y)\neq 0$ for every $y\in\Gamma$. If we set $X:=\R^d\backslash\Gamma^\circ$ it follows from Theorem \ref{characterization of p-convexity in the complement of cones} i) that $X$ is $P$-convex for supports. Since $P$ is hypoelliptic $X$ is $P$-convex for singular supports as well. But because $x\in\Gamma$ and $\sigma^0_P(x)=0$ by inequality (\ref{key inequality}) $X\times\R$ is not $P^+$-convex for singular supports by Theorem \ref{characterization of p-convexity in the complement of cones}. Thus, for the hypoelliptic polynomial $P$ we have
\[P(D):\mathscr{D}'(X)\rightarrow \mathscr{D}'(X)\mbox{ is surjective}\]
but
\[P^+(D):\mathscr{D}'(X\times\R)\rightarrow \mathscr{D}'(X\times\R)\mbox{ is not surjective.}\]
Since for hypoelliptic $P$ the kernels of
\[P(D):C^\infty(X)\rightarrow C^\infty(X)\mbox{ and }P(D):\mathscr{D}'(X)\rightarrow\mathscr{D}'(X)\]
coincide as locally convex spaces it is a Fr\'echet-Schwartz space. Hence it has $(\Omega)$ if and only if it has $(P\Omega)$. By the results of Bonet and Doma\'nski mentioned in the introduction the above polynomial $P$ and the open set $X$ gives thus a surjective hypoelliptic differential operator
\[P(D):C^\infty(X)\rightarrow C^\infty(X)\]
such that its kernel does not have property $(\Omega)$. This should be compared with Vogt's classical result \cite{Vogt} that the kernel of an elliptic differential operators always has $(\Omega)$.

\begin{small}
{\sc FB 4 - Mathematik, Universit\"at Trier,
54286 Trier, Germany}

{\it E-mail address: kalmesth@uni-trier.de}
\end{small}

\end{document}